\documentstyle[12pt]{article}
\let\text\mbox
\def\limfunc#1{\mbox{\rm #1}\,}
\begin{document}
\author{{\bf Steven Duplij}\thanks{Alexander von Humboldt Fellow}
\thanks{On leave of absence from: {\sl Theory Division,
Nuclear Physics Laboratory,
Kharkov State University, KHARKOV 310077, Ukraine}}
\thanks{E-mail: duplij@physik.uni-kl.de} \\
{\sl Department of Physics, University of Kaiserslautern,}\\
{\sl Postfach 3049, D-67653 Kaiserslautern, Germany}}
\title{{\bf SECOND $N=1$\ SUPERANALOG OF COMPLEX STRUCTURE}}
\date{October 9, 1995}
\maketitle

\begin{abstract}
We found another $N=1$ odd superanalog of complex structure (the even one is
widely used in the theory of super Riemann surfaces). New $N=1$
superconformal-like transformations are similar to anti-holomorphic ones of
nonsupersymmetric complex function theory. They are dual to the ordinary
superconformal transformations subject to the Berezinian addition formula
presented, noninvertible, highly degenerated and twist parity of the tangent
space in the standard basis. They also lead to the ''mixed cocycle
condition'' which can be used in building noninvertible objects analogous to
super Riemann surfaces. A new parametrization for the superconformal group
is presented which allows us to extend it to a semigroup and to unify the
description of old and new transformations. \vskip10pt \noindent PACS
numbers: 11.30.Pb; 02.20Mp
\end{abstract}
\begin{flushright}KL-TH-95/23
\end{flushright}
\newpage\ The idea of superconformal symmetry is exceptionally important in
the theory of super Riemann surfaces \cite{dho/pho1} and in two-dimensional
superconformal field theories \cite{kaku1}. The main and fundamental
ingredient of the idea is a special class of reduced mappings of
two-dimensional $(1|1)$ complex superspace, namely superconformal
transformations \cite{cra/rab}. In the local approach to super Riemann
surfaces represented as collections of open superdomains the superconformal
transformations are used as gluing transition functions \cite{cra/rab,fri}.
{}From another side they appear as a result of the special reduction of the
structure supergroup \cite{gid/nel1}. Here we consider an alternative
tangent space reduction, which leads to new transformations (see also \cite
{:2a}\nocite{dup6,dup7} and \cite{dup12}).

We use the functional approach to superspace \cite{rog1} which admits
existence of nontrivial topology in odd directions \cite{rab/cra1} and can
be suitable for physical applications \cite{:2b}\nocite{bru,van}. Also we
exploit the coordinate language which is more physically transparent
and adequate in constructing objects having new features.

Locally $\left( 1|1\right) $-dimensional superspace ${\rm C}^{1|1}$ is
described by $Z=\left( z,\theta \right) $, where $z$ is an even coordinate
and $\theta $ is an odd one. The most intriguing peculiarity of the
functional definition of superspace \cite{rog1} is the existence of soul
parts in the even coordinate $z=z_{body}+z_{soul},\;z_{body}=\epsilon \left(
z\right) ,\;z_{soul}\stackrel{def}{=}z-z_{body}$, where $\epsilon $ is a
body map \cite{rog1} vanishing all nilpotent generators. The body map acts
on the coordinates as follows $\epsilon \left( z\right) =z_{body},\;\epsilon
\left( \theta \right) =0$. This allows one to consider non-trivial soul
topology in even directions on a par with odd ones \cite{rab/cra1}.

A superanalytic (SA) transformation ${\rm T}_{SA}:{\rm C}^{1|1}\rightarrow
{\rm C}^{1|1}$ is
\begin{equation}
\label{1}
\begin{array}{ccc}
\tilde z & = & f\left( z\right) +\theta \cdot \chi \left( z\right) , \\
\tilde \theta  & = & \psi \left( z\right) +\theta \cdot g\left( z\right) ,
\end{array}
\end{equation}

\noindent where four component functions $f\left( z\right) ,g\left( z\right)
:{\rm C}^{1|0}\rightarrow {\rm C}^{1|0}$ and $\psi \left( z\right) ,\chi
\left( z\right) :{\rm C}^{1|0}\rightarrow {\rm C}^{0|1}$ satisfy some
supersmooth conditions generalizing $C^\infty $ \cite{rog1}, and
simultaneously they can be noninvertible \cite{:2a} (here and in the
following we denote even functions and variables by Latin letters and odd
ones by Greek letters, point is a product in Grassmann algebra). The set of
invertible and noninvertible SA transformations (\ref{1}) form a semigroup
of superanalytic transformations ${\cal T}_{SA}$ \cite{:2a}. The invertible
transformations are in its subgroup, while the noninvertible ones are in an
ideal (see \cite{:2a} and \cite{:3a}\nocite{dup11,dup8,dup5} for details).

The invertibility of the superanalytic transformation (\ref{1}) is
determined first of all by invertibility of the even functions $f\left(
z\right) $ and $g\left( z\right) $, because odd functions are noninvertible
by definition. In case $\epsilon \left( g\left( z\right) \right) \neq 0$ for
SA transformations (\ref{1}) the superanalog of a Jacobian, the Berezinian
\cite{berezin}, can be determined $\limfunc{Ber}\left( \tilde Z/Z\right) =%
\frac{f^{\prime }\left( z\right) }{g\left( z\right) }+\frac{\chi \left(
z\right) \psi ^{\prime }\left( z\right) }{g^2\left( z\right) }+\theta \left(
\frac{\chi \left( z\right) }{g\left( z\right) }\right) ^{\prime }$, where
prime is a differentiation by argument (or by $z$). Therefore, we can
classify the transformations (\ref{1}) in the following way: 1) the
Berezinian exists and invertible ($\epsilon \left( g\left( z\right) \right)
\neq 0,\;\epsilon \left( f\left( z\right) \right) \neq 0$); 2) the
Berezinian exists and noninvertible ($\epsilon \left( g\left( z\right)
\right) \neq 0,\;\epsilon \left( f\left( z\right) \right) =0$); 3) the
Berezinian does not exist ($\epsilon \left( g\left( z\right) \right)
=0,\;\epsilon \left( f\left( z\right) \right) =0$). First type of SA
transformations form a subgroup of the superanalytic semigroup, while second
two types are in an ideal of the semigroup \cite{:2a}.

The tangent superspace in ${\rm C}^{1|1}$ is defined by the standard basis $%
\left\{ \partial ,\;D\right\} $, where $D=\partial _\theta +\theta \partial $%
, $\partial _\theta =$ $\partial /\partial \theta ,\;\partial =$ $\partial
/\partial z$. The dual cotangent space is spanned by 1-forms $\left\{
dZ,\;d\theta \right\} $, where $dZ=dz+\theta d\theta $ (the signs as in \cite
{cra/rab}). In these notations the supersymmetry relations are $D^2=\partial
,\;dZ^2=dz$. The semigroup of SA transformations acts in the tangent and
cotangent superspaces by means of the tangent space matrix $P_A$ as $\left(
\begin{array}{c}
\partial \\
D
\end{array}
\right) =P_A\left(
\begin{array}{c}
\tilde \partial \\
\tilde D
\end{array}
\right) $ and $\left(
\begin{array}{cc}
d\tilde Z, & d\tilde \theta
\end{array}
\right) =\left(
\begin{array}{cc}
dZ, & d\theta
\end{array}
\right) P_A$, where

\begin{equation}
\label{m}P_A=\left(
\begin{array}{cc}
\partial \tilde z-\partial \tilde \theta \cdot \tilde \theta & \partial
\tilde \theta \\
D\tilde z-D\tilde \theta \cdot \tilde \theta & D\tilde \theta
\end{array}
\right) .
\end{equation}

In case of invertible SA transformations the matrix $P_A$ defines structure
of a supermanifold for which these transformations play the part of
transition functions, and $\limfunc{Ber}\left( \tilde Z/Z\right) =\limfunc{%
Ber}P_A.$ Therefore different reductions of the matrix $P_A$ give us various
additional supermanifold structures \cite{gid/nel1}. It was shown in \cite
{dup12} that there exist two nontrivial reductions of any supermatrix $P_A$.
Indeed, if $\epsilon \left( D\tilde \theta \right) \neq 0$ we observe that $%
\limfunc{Ber}P_A=\frac{\partial \tilde z-\partial \tilde \theta \cdot \tilde
\theta }{D\tilde \theta }+\frac{\left( D\tilde z-D\tilde \theta \cdot \tilde
\theta \right) \partial \tilde \theta }{\left( D\tilde \theta \right) ^2}$.
Then using the Berezinian addition theorem \cite{dup12} we obtain the formula

\begin{equation}
\label{b}\limfunc{Ber} P_A=\limfunc{Ber} P_S+\limfunc{Ber} P_T,
\end{equation}

\noindent where
\begin{equation}
\label{ms}P_S\stackrel{def}{=}\left(
\begin{array}{cc}
\partial \tilde z-\partial \tilde \theta \cdot \tilde \theta & \partial
\tilde \theta \\
0 & D\tilde \theta
\end{array}
\right) ,
\end{equation}
\begin{equation}
\label{mt}P_T\stackrel{def}{=}\left(
\begin{array}{cc}
0 & \partial \tilde \theta \\
D\tilde z-D\tilde \theta \cdot \tilde \theta & D\tilde \theta
\end{array}
\right) .
\end{equation}

Denote sets of the matrices (\ref{ms}) and (\ref{mt}) by ${\bf P}_S$ and $%
{\bf P}_T$ respectively. Then their intersection ${\bf P}_D={\bf P}_S\cap
{\bf P}_T$ is a set of the degenerated matrices $P_D$ of the form
\begin{equation}
\label{md}P_D\stackrel{def}{=}\left(
\begin{array}{cc}
0 & \partial \tilde \theta \\
0 & D\tilde \theta
\end{array}
\right) ,
\end{equation}

\noindent which depend on the odd coordinate $\theta $ transformation only.
The degenerated matrix of the shape (\ref{md}) can be obtained by projection
from $P_S$ and $P_T$ matrices using the following equations

\begin{equation}
\label{tpt}Q\stackrel{def}{=}\partial \tilde z-\partial \tilde \theta \cdot
\tilde \theta =0,
\end{equation}

\begin{equation}
\label{scf}\Delta \stackrel{def}{=}D\tilde z-D\tilde \theta \cdot \tilde
\theta =0
\end{equation}

\noindent correspondingly. It means that, if the transformation of the odd
sector (second line in (\ref{1})) is given, i.e. the functions $\psi \left(
z\right) $ and $g\left( z\right) $ are fixed, the conditions (\ref{tpt}) and
(\ref{scf}) determine behavior of the even sector (functions $f\left(
z\right) $ and $\chi \left( z\right) $). In this case, since the degenerated
matrix $P_D$ depends on the odd sector transformation only, we obtain
\begin{equation}
\label{deg}P_D=P_S|_{Q=0}=P_T|_{\Delta =0}.
\end{equation}

An opposite situation occurs if we apply the conditions (\ref{tpt}) and (\ref
{scf}) to the matrices $P_S$ and $P_T$ in a reverse order. Then we derive
\begin{equation}
\label{4}P_{SCf}\stackrel{def}{=}P_S|_{\Delta =0}
\end{equation}
\begin{equation}
\label{5}P_{TPt}\stackrel{def}{=}P_T|_{Q=0}.
\end{equation}

The condition $\Delta =0$ (\ref{scf}) gives us superconformal (SCf )
transformations ${\rm T}_{SCf}$ \cite{cra/rab} and the reduced matrix $%
P_{SCf}$ (\ref{4}) is a result of the standard reduction of structure
supergroup (in the invertible case \cite{gid/nel1}). Another condition $%
\Delta =0$ (\ref{tpt}) leads to the degenerated transformations ${\rm T}%
_{TPt}$ twisting parity of the standard tangent space (TPt ) \cite{:2a}.
The alternative reduction \cite{dup12} of the tangent space supermatrix $P_A$
gives us the supermatrix $P_{TPt}$ (\ref{5}). The dual role of SCf and TPt
transformations is clearly seen from the Berezinian addition theorem (\ref{b}%
) (see \cite{dup12}) and the projections (\ref{4}) and (\ref{5}). Since SCf
transformations give us a superanalog of complex structure \cite{levn1,schw3}%
, we can treat TPt transformations as another odd $N=1$ superanalog of
complex structure in a certain extent\footnote{%
It is more natural to call TPt transformations anti-SCf transformations due
to the following analogy with the nonsupersymmetric case. For an ordinary $%
2\times 2$ matrix $P=\left(
\begin{array}{cc}
a & b \\
c & d
\end{array}
\right) $ we obviously have the following identity $\det P=\det \left(
\begin{array}{cc}
a & 0 \\
0 & d
\end{array}
\right) +\det \left(
\begin{array}{cc}
0 & b \\
c & 0
\end{array}
\right) =\det P_{Diag}+\det P_{Antidiag}$, which can be called a
''determinant addition formula''. In the complex function theory the first
matrix describes the tangent space matrix of holomorphic mappings and the
second one---of antiholomorphic mappings. In supersymmetric case the
supermatrices $P_S$ and $P_T$ play the role similar to one of the
nonsupersymmetric diagonal and antidiagonal matrices in ordinary theory as
it is seen from (\ref{b}). Therefore, if $P_{SCf}$ generalizes the tangent
space matrix of holomorphic mappings, supermatrices $P_{TPt}$ could be
considered as respective generalization for antiholomorphic mappings.}.

Using (\ref{4}) and (\ref{5}) with the obvious relation $\limfunc{Ber}P_D=0$
we can project the Berezinian addition equality (\ref{b}) to ${\rm T}_{SCf}$
and ${\rm T}_{TPt}$ as follows

\begin{equation}
\label{bu}\limfunc{Ber}P_A=\left\{
\begin{array}{cc}
\limfunc{Ber}P_{SCf}, & \Delta =0, \\
\limfunc{Ber}P_{TPt}, & Q=0.
\end{array}
\right.
\end{equation}

A general relation between $Q$ and $\Delta $ is $Q-D\Delta =\left( D\tilde
\theta \right) ^2$. After corresponding projections we have
\begin{equation}
\label{q}Q|_{\Delta =0}=\left( D\tilde \theta \right) ^2,\text{ (SCf ),}
\end{equation}
\begin{equation}
\label{dl}\Delta |_{Q=0}\equiv \Delta _0=\partial _\theta \tilde z-\partial
_\theta \tilde \theta \cdot \tilde \theta ,\text{ (TPt )}.
\end{equation}

It is remarkable to notice the similarity of (\ref{tpt}) and (\ref{dl}).
Using (\ref{q}) one obtains \cite{gid/nel1}

\begin{equation}
\label{Pscf}P_{SCf}=\left(
\begin{array}{cc}
\left( D\tilde \theta \right) ^2 & \partial \tilde \theta \\
0 & D\tilde \theta
\end{array}
\right) .
\end{equation}

If $\varepsilon \left( D\tilde \theta \right) \neq 0$ then $\limfunc{Ber}%
P_{SCf}$ can be determined and it is

\begin{equation}
\label{berPscf}\limfunc{Ber}P_{SCf}=D\tilde \theta .
\end{equation}

In case $\varepsilon \left( D\tilde \theta \right) =0$ the Berezinian cannot
be defined, but we can accept (\ref{berPscf}) as a definition of the
Jacobian of noninvertible SCf transformations (see \cite{:2a} and
\cite{:2d}\nocite{dup4,dup10}).

{}From (\ref{dl}) we derive

\begin{equation}
\label{Ptpt}P_{TPt}=\left(
\begin{array}{cc}
0 & \partial \tilde \theta \\
\partial _\theta \tilde z-\partial _\theta \tilde \theta \cdot \tilde \theta
& D\tilde \theta
\end{array}
\right)
\end{equation}

\noindent (cf. (\ref{ms})). If $\varepsilon \left( D\tilde \theta \right)
\neq 0$ the Berezinian of $P_{TPt}$ can be determined as

\begin{equation}
\label{berPtpt}\limfunc{Ber}P_{TPt}=\frac{\Delta _0\cdot \partial \tilde
\theta }{\left( D\tilde \theta \right) ^2}.
\end{equation}

{}From (\ref{dl}) it follows that $D\Delta _0=-\left( D\tilde \theta \right)
^2 $ and, therefore, $\partial \Delta _0=-2\cdot D\tilde \theta \cdot
\partial \tilde \theta $, which gives $\limfunc{Ber}P_{TPt}=\frac{\partial
\Delta _0\cdot \Delta _0}{2\left( D\tilde \theta \right) ^3}$. Since $\Delta
_0$ is odd and so nilpotent, $\limfunc{Ber}P_{TPt}$ is also nilpotent and
pure soul. The Berezinian can be also presented as
\begin{equation}
\label{berPtpt2}\limfunc{Ber}P_{TPt}=D\left( \frac{D\tilde z}{D\tilde \theta
}\right)
\end{equation}

\noindent which should be remarkably compared with (\ref{berPscf}).

The most intriguing peculiarity of TPt transformations is twisting the
parity of tangent and cotangent spaces in the standard basis, viz.
\begin{equation}
\label{8}\text{SCf: }\left\{
\begin{array}{ccc}
D & = & \left( D\tilde \theta \right) \cdot \tilde D, \\
d\tilde Z & = & \left( D\tilde \theta \right) ^2\cdot dZ,
\end{array}
\right. \;\text{\ \ \ \ TPt: }\left\{
\begin{array}{ccc}
\partial & = & \partial \tilde \theta \cdot \tilde D, \\
d\tilde Z & = & \Delta _0\cdot d\theta .
\end{array}
\right.
\end{equation}

The reduction conditions (\ref{tpt}) and (\ref{scf}) fix 2 of 4 component
functions form (\ref{1}) in each case. Usually \cite{cra/rab} SCf
transformations ${\rm T}_{SCf}$ are parametrized by $\left(
\begin{array}{c}
f \\
\psi
\end{array}
\right) $, while other functions are found from (\ref{tpt}) and (\ref{scf}).
However, the latter can be done for invertible transformations only. To
avoid this difficulty we introduce an alternative parametrization by the
pair $\left(
\begin{array}{c}
g \\
\psi
\end{array}
\right) $, which allows us to consider SCf and TPt transformations in a
unified way and include noninvertibility. Indeed, fixing $g\left( z\right) $
and $\psi \left( z\right) $ we find for other component functions of (\ref{1}%
) the equations
\begin{equation}
\label{9}\left\{
\begin{array}{ccc}
f_n^{\prime }\left( z\right) & = & \psi ^{\prime }\left( z\right) \psi
\left( z\right) +
\frac{1+n}2g^2\left( z\right) , \\ \chi _n^{\prime }\left( z\right) & = &
g^{\prime }\left( z\right) \psi \left( z\right) +ng\left( z\right) \psi
^{\prime }\left( z\right) ,
\end{array}
\right.
\end{equation}

\noindent where $n=\left\{
\begin{array}{cc}
+1, & \text{SCf,} \\ -1, & \text{TPt,}
\end{array}
\right. $ can be treated as a projection of some ''reduction spin''
switching the type of transformation. So the reduced transformation of the
even coordinate (see (\ref{1})) should contain this additional index, i.e. $%
z\rightarrow \tilde z_n$ (at this point the analogy with complex structure is
most transparent).
Since $f_{-1}^{\prime }\left( z\right) =\psi
^{\prime }\left( z\right) \psi \left( z\right) $ is nilpotent, TPt
transformations are always noninvertible and high degenerated after the body
mapping. The unified multiplication law is
\begin{equation}
\label{10}\left(
\begin{array}{c}
h \\
\varphi
\end{array}
\right) _n*\left(
\begin{array}{c}
g \\
\psi
\end{array}
\right) _m=\left(
\begin{array}{c}
g\cdot h\circ f_m+\chi _m\cdot \psi \cdot h^{\prime }\circ f_m+\chi _m\cdot
\varphi ^{\prime }\circ f_m \\
\varphi \circ f_m+\psi \cdot h\circ f_m
\end{array}
\right) ,
\end{equation}

\noindent where ($*$) is transformation composition and ($\circ $) is
function composition. For ''reduction spin'' projections we have only two
definite products $\left( +1\right) *\left( +1\right) =\left( +1\right) $
and $\left( +1\right) *\left( -1\right) =\left( -1\right) $. The first
formula is a consequence of ${\bf P}_S\cdot {\bf P}_S\subseteq {\bf P}_S$
(see (\ref{ms})), which is simple manifestation of the fact that SCf
transformations ${\rm T}_{SCf}$ form a substructure \cite{gid/nel1}, i.e. a
subsemigroup ${\cal T}_{SCf}$ of SA\ semigroup ${\cal T}_{SA}$ (in the
invertible case---a subgroup \cite{cra/rab}). From ${\bf P}_S\cdot {\bf P}%
_S\subseteq {\bf P}_S$ it also follows the standard (for component functions
too) cocycle condition \cite{cra/rab} ${\rm \tilde T}_{SCf}*{\rm T}_{SCf}=%
\widetilde{{\rm \tilde T}}_{SCf}$ (having identical arrows, i.e. (SCf)
actions) on triple overlaps $U\cap \tilde U\cap \widetilde{\tilde U}$, where
$U$, $\tilde U$, $\widetilde{\tilde U}$ are open superdomains and ${\rm T}%
:U\rightarrow \tilde U,\;{\rm \tilde T}:\tilde U\rightarrow \widetilde{%
\tilde U},\;\widetilde{{\rm \tilde T}}:U\rightarrow \widetilde{\tilde U}$.
In the invertible SCf case the cocycle condition leads to the definition of
a super Riemann surface as a holomorphic $\left( 1|1\right) $-dimensional
supermanifold equipped with an additional one-dimensional subbundle \cite
{cra/rab,gid/nel1,levn1}, which grounds on the cocycle relation $D\widetilde{%
\tilde \theta }=D\tilde \theta \cdot \tilde D\widetilde{\tilde \theta }$ and
the formula (\ref{berPscf}). Unfortunately, TPt transformations ${\rm T}%
_{TPt}$ form a subsemigroup only providing additional conditions on
component functions \cite{:2a}. However, they have also another important
abstract meaning: using the unrestricted relation ${\bf P}_T\cdot {\bf P}%
_S\subseteq {\bf P}_T$ we obtain a ''mixed cocycle condition'' ${\rm \tilde T%
}_{SCf}*{\rm T}_{TPt}=\widetilde{{\rm \tilde T}}_{TPt}$ (having different
arrows). Then we derive the ''mixed cocycle relation''
\begin{equation}
\label{mx}\partial \widetilde{\tilde \theta }=\partial \tilde \theta \cdot
\tilde D\widetilde{\tilde \theta }
\end{equation}
which can be exploited in constructing new objects analogous to super
Riemann surfaces. It is remarkable that under the degenerated (Deg)
transformations defined by (\ref{deg}) the both cocycle relations hold valid
simultaneously. Also, Deg transformations form a subsemigroup ${\cal T}_{Deg}
$ in ${\cal T}_{SA}$, because of ${\bf P}_D\cdot {\bf P}_D\subseteq {\bf P}_D
$. Moreover, ${\cal T}_{Deg}$ is an ideal in ${\cal T}_{SA}$, ${\cal T}_{SCf}
$ and ${\cal T}_{TPt}$ since ${\bf P}_D\cdot {\bf P}_A\subseteq {\bf P}_D$, $%
{\bf P}_D\cdot {\bf P}_S\subseteq {\bf P}_D$ and ${\bf P}_D\cdot {\bf P}%
_T\subseteq {\bf P}_D$. The degenerated transformations are characterized by
one odd function $\psi \left( z\right) $ only and by the absence of the $%
\theta $-dependence of the transformation $Z\rightarrow \tilde Z$ (see (\ref
{dl})), so that $\tilde z_{Deg}=f\left( z\right) $, $\tilde \theta
_{Deg}=\psi \left( z\right) $ and $f^{\prime }\left( z\right) =\psi ^{\prime
}\left( z\right) \psi \left( z\right) $. The multiplication in ${\cal T}%
_{Deg}$ coincides with the second row of (\ref{10}).

We conclude that thorough consideration of invertibility, while
super-generalizing standard constructions of string theory,
leads to  nontrivial consequences and possibility of building
new objects analogous, for instance, to super Riemann surfaces, which could
give additional contributions to fermionic string amplitude. It would be
also interesting to work out sequences of
noninvertible functions and corresponding bundles or their generalizations.

Author would like to thank M. Grisaru, P. Howe, J. Kupsch, P. van
Nieuwenhuizen, W. R\"uhl and P. Townsend for useful discussions.

\end{document}